\documentclass[review]{elsarticle}
\usepackage{algorithm, algorithmic}
\usepackage{hyperref}
\usepackage{lineno}
\modulolinenumbers[5]
\usepackage{amsmath}
\journal{Journal of \LaTeX\ Templates}

\begin{document}

\begin{frontmatter}

\title{A Hybrid Numerical Algorithm for Evaluating n-th Order Tridiagonal Determinants}

\author{Moawwad El-Mikkawy\corref{mycorrespondingauthor}}
\author{Abdelrahman Karawia}
\address{Mathematics Department, Faculty of Science, Mansoura University, Mansoura, 35516, Egypt}
\cortext[mycorrespondingauthor]{Corresponding author}




\begin{abstract}
The principal minors of a tridiagonal matrix satisfy two-term and three-term recurrences \cite{Elmikkawy2004,Elmikkawy2003}. Based on these facts, the current article presents a new efficient and reliable hybrid numerical algorithm for evaluating general n-th order tridiagonal determinants in linear time. The hybrid numerical algorithm avoid all symbolic computations. The algorithm is suited for implementation using computer languages such as FORTRAN, PASCAL, ALGOL, MAPLE, MACSYMA and MATHEMATICA. Some illustrative examples are given. Test results indicate the superiority of the hybrid numerical algorithm.
\end{abstract}


\end{frontmatter}

\linenumbers

\section{Introduction}
A general tridiagonal matrix $T_n=(t_{ij})_{1\le i,j\le n}$ takes the form:
\begin{equation}
    T_n=(t_{ij})= \left[
\begin{array}{cccccc}
    d_1 & a_1 & 0 & \cdots & \cdots & 0\\
   b_1 & d_2 & a_2 &\ddots &  & \vdots \\
   0 & b_2 & d_3 & \ddots & 0 & \vdots\\
   \vdots & \ddots &\ddots &\ddots &\ddots & 0\\
   \vdots & & 0 &\ddots &\ddots & a_{n-1}\\ 
   0& \cdots & \cdots & 0 & b_{n-1} & d_n\\
 \end{array}
\right]_n, n\ge 3
\end{equation}
in which $t_{ij}=0$ whenever $|i-j| > 1$.\\
These matrices arise frequently in a wide range of scientific and engineering fields \cite{MOLINARI2008,Jina20013,JIA2015,Feng2021}. For instance, telecommunication, parallel computing, and statistics. For the matrix $T_n$ in (1), there is no need to store the zero elements. Consequently, we can use three vectors $\mathbf{a}=(a_1, a_2, \cdots, a_{n-1})$, $\mathbf{b}=(b_1, b_2, \cdots, b_{n-1})$, and $\mathbf{d}=(d_1, d_2, \cdots, d_n)$ to store the non-zero elements of $T_n$ in $3n-2$ memory locations rather than $n^2$ for a full matrix. This is always a good habit in computation. When we consider the matrix $T_n$ in (1), it is useful to add an additional n-dimensional vector $\mathbf{c}=(c_1, c_2, \cdots, c_n)$  given by:
\begin{equation}
    c_i= \left\{
  \begin{array}{lr} 
      d_1 & \text{if}\quad i=1, \\
      d_i-a_{i-1} b_{i-1}/c_{i-1} & \text{if}\quad i=2, 3, \cdots, n 
      \end{array}
      \right.
\end{equation}
By adding this vector $\mathbf{c}$ , we are able to: \\
(i) evaluate $det(T_n)$ in linear time \cite{Elmikkawy2004},\\
(ii) write down the Doolittle and Crout LU factorizations of the matrix $T_n$ \cite{Burden2001}, and \\
(iii) check whether or not a symmetric tridiagonal matrix $T_n$ is the positive definite. In fact if $T_n$ is symmetric then it is positive definite if and only if $c_i>0, i=1, 2, \cdots, n$ \cite{Burden2001}.\\
\\
In \cite{url2018:online}, the following question has been raised:\\
Is there a fast way to prove that the tridiagonal matrix\\
$$A=\left[
\begin{array}{ccccc}
    4 & 2 & 0 & 0&0\\
   2 & 5 & 2 &0 & 0 \\
   0 & 2 & 5 & 2&0\\
   0& 0 & 2 & 5 & 2\\
   0& 0 & 0 & 2 & 5\\
 \end{array}
\right]$$\\
is a positive definite?\\
Our answer is: $A$ is actually positive definite since $c_i=4 > 0, i=1, 2, 3, 4, 5$ as can be easily checked. This is the easiest way to check the positive definiteness of a symmetric tridiagonal matrix.\\

The current article is organized as follows. The main result is presented in Section 2. In Section 3, numerical tests and some illustrative examples are given. The conclusion is presented in Section 4.
\section{The Main Result}
This section is mainly devoted to constructs a hybrid numerical algorithm for evaluating n-th order tridiagonal determinant of the form (1). \\
Let:
\begin{equation}
    f_1=|d_1|=d_1, f_i=\left|
\begin{array}{cccccc}
    d_1 & a_1 & 0 & \cdots & \cdots & 0\\
   b_1 & d_2 & a_2 &\ddots &  & \vdots \\
   0 & b_2 & d_3 & \ddots & 0 & \vdots\\
   \vdots & \ddots &\ddots &\ddots &\ddots & 0\\
   \vdots & & 0 &\ddots &\ddots & a_{i-1}\\ 
   0& \cdots & \cdots & 0 & b_{i-1} & d_i\\
 \end{array}
\right|, i=2, 3, \cdots. n
\end{equation}
Therefore, $f_i, i=1, 2, \cdots, n$ are the principal minors of $T_n$.
The determinants in (3) satisfy a two-term recurrence \cite{Elmikkawy2003}
\begin{equation}
    f_i=\prod_{r=1}^ic_r=c_i f_{i-1}, i=1, 2, \cdots, n, f_0=1
\end{equation}
The three-term recurrence
\begin{equation}
    f_i=d_i f_{i-1}-a_{i-1}b_{i-1}f_{i-2}, i=2, 3, \cdots, n, f_0=1, f_1=d_1
\end{equation}
is also valid \cite{Elmikkawy2003}.\\
For convenience of the reader it is convenient to describe the \textbf{DETGTRI} algorithm in which $z$ is just a symbolic name \cite{Elmikkawy2004}.
\begin{algorithm}[!h]
  \caption{\textbf{DETGTRI}}
  \textbf{Input:} $n$ and the components of the vectors $\mathbf{a}, \mathbf{b}$, and $\mathbf{d}$.\\
  \textbf{Output:} $det(T_n)$.\\
  \textbf{Step 1:} For $k$ from $2$ to $n$ do\\
  \hspace*{2cm}  Compute and simplify:\\ 
  \hspace*{2cm} If $d_{k-1}=0$ then $d_{k-1}=z$ end if.\\
  \hspace*{2.1cm}$d_k := d_k - a_{k-1} b_{k-1}/d_{k-1}$\\
  \hspace*{1.3cm} End do\\
   \textbf{Step 2:} Compute  $P(z) =\prod_{r=1}^n d_r$\\
   \textbf{Step 3:} Set $det(T_n)=P(0)$.
  \label{alg1}
  \end{algorithm}
\\
\noindent Based on the three-term recurrence (5), we may formulate the following algorithm.\\
\begin{algorithm}[!h]
  \caption{}
  \textbf{Input:} $n$ and the components of the vectors $\mathbf{a}, \mathbf{b}$, and $\mathbf{d}$.\\
  \textbf{Output:} $det(T_n)$.\\
  \textbf{Step 1:} Set $f_0=1$ and $f_1=d_1$.\\
  \textbf{Step 2:} For $i$ from $2$ to $n$ do\\
  \hspace*{2cm}  $f_i=d_if_{i-1} - a_{i-1}b_{i-1}f_{i-2}$,\\
  \hspace*{1.5cm}End do.\\
  \textbf{Step 3:} Set $det(T_n)=f_n$.\\
  \label{alg2}
  \end{algorithm}
\\
\\
\noindent At this stage, we present the following hybrid numerical algorithm.\\
\begin{algorithm}[!h]
  \caption{}
  \textbf{Input:} $n$ and the components of the vectors $\mathbf{a}, \mathbf{b}$, and $\mathbf{d}$.\\
  \textbf{Output:} $det(T_n)$.\\
  \textbf{Step 1:} Set $c_1 = d_1$, $f_1 = d_1$, and $m = 1$.\\
  \textbf{Step 2:} While $m\le n-1$ and $c_m \ne 0$ do\\
 \hspace*{2cm} $m = m+1$,\\
  \hspace*{2cm}  $c_m = d_m-a_{m-1}b_{m-1}/c_{m-1}$,\\
   \hspace*{2cm}$f_m = c_mf_{m-1}$,\\
  \hspace*{1.5cm}End do.\\
  \textbf{Step 3:} For $k = m+1$ to $n$ do\\
  \hspace*{2cm} $f_k = d_kf_{k-1} - a_{k-1}b_{k-1}f_{k-2}$\\
  \hspace*{1.5cm}End do.\\
  \textbf{Step 4:}
  Set $det(T_n)=f_n$.\\
  \label{alg3}
  \end{algorithm}
  \\

The hybrid numerical algorithm has the same computational cost as the algorithms \textbf{DETGTRI} and \textbf{Algorithm 2}. \textbf{Algorithm 3} links two methods and has the advantage that no symbolic computations are involved.\\
\textbf{Remark:} It should be noted that \textbf{Step 3} in \textbf{Algorithm 3} is redundant and will not be executed at all if $c_i\ne 0, i = 1, 2, \cdots, n-1$. Therefore, we only need \textbf{Step 1}, \textbf{Step 2} and \textbf{Step 4}. For positive definite and strictly diagonally dominant matrices,
this is always the case. The implementation of the hybrid numerical algorithm using any computer language are straight forward.
\section{Numerical Tests and Illustrative Examples}
In this section, we are going to consider Some numerical tests and illustrative examples. All computations are carried out using laptop machine with a 2.50GHz CPU, 8GB of RAM, AMD A10-9620P RADEON R5 processor and Maple 2021.\\
\textbf{Example 3.1.} Consider the tridiagonal matrix $T_n$, with $n=4$ given by:\\
\\
\hspace*{3cm} $T_n=(t_{ij})=\left[
\begin{array}{cccc}
    1 & 1 & 0 & 0\\
   1 & 1 & -1 &0 \\
   0 & 1 & 2 & 1\\
   0& 0 & -3& -1\\
 \end{array}
\right]_4$\\
Find $det(T_n)$.\\
\textbf{Solution}:\\
We have:\\
$a_1=1, a_2=-1, a_3=1, b_1=1, b_2=1, b_3=-3, d_1=1, d_2=1, d_3=2,$ and $d_4=-1$.\\
By applying the \textbf{Algorithm 3}, we obtain\\
\textbf{Step 1:} $c_1=d_1=1$, $f_1=d_1=1$ and $m=1$\\
\textbf{Step 2:} $m=2, c_2=0, f_2=c_2 f_1=0$.\\
\textbf{Step 3:} $f_3=d_3f_2-a_2b_2f_1=(2)(0)-(-1)(1)(1)=1,$ and $f_4=d_4f_3-a_3b_3f_2=(-1)(1)-(1)(-3)(0)=-1$.\\
\textbf{Step 4:} $det(T_n)=f_4=-1$.\\
\\
\textbf{Example 3.2.} Consider $T_n,$ with $n=9$, given by:\\
$a_i=-1, b_i=-1, d_i=2, i=1, 2, \cdots, n-1,$ and $d_{n}=2$.\\
By applying the \textbf{Algorithm 3}, we get\\
\textbf{Step 1:} $c_1=2$, and $f_1=2$.\\
\textbf{Step 2:} \\
$$\begin{tabular}{ccccccccc}
\hline
$m$ & 2 & 3 & 4 & 5 & 6 & 7 & 8 & 9 \\
\hline
$c_m$ & $\frac{3}{2}$ & $\frac{4}{3}$ & $\frac{5}{4}$  & $\frac{6}{5}$ & $\frac{7}{6}$ & $\frac{8}{7}$ &$\frac{9}{8}$ &$\frac{10}{9}$\\
$f_m$ & 3 & 4 & 5 & 6 & 7 & 8 & 9 & 10\\
\hline
\end{tabular}$$\\
\textbf{Step 4:} $det(T_n)=f_9=10$.\\
\\
\textbf{Example 3.3.} Let $T_n$ is given by:\\
\hspace*{3cm} $T_n=(t_{ij})\left[
\begin{array}{cccccc}
    1 & 1 & 0 & \cdots & \cdots & 0\\
   1 & 1 & 1 &\ddots &  & \vdots \\
   0 & 1 & 1 & \ddots & 0 & \vdots\\
   \vdots & \ddots &\ddots &\ddots &\ddots & 0\\
   \vdots & & 0 &1 &1 & 1\\ 
   0& \cdots & \cdots & 0 & 1 & 1\\
 \end{array}
\right]_n$
\\
By using (4), we get:\\
$det(T_n)= \left\{
  \begin{array}{lr} 
      1 & \text{if}\quad n\equiv 0\quad or\quad 1\quad mod(6), \\
      0 & \text{if}\quad n\equiv 2\quad or\quad 5\quad mod(6), \\
      -1 & \text{if}\quad n\equiv 3\quad or\quad 4\quad mod(6).
      \end{array}
      \right.
$\\
\\
Now, it is time to consider Example 3.3 as a test problem to compare the three algorithms and the MATLAB function $det()$. For these algorithms, we get the results presented in Table 1. The \textbf{DETGTRI} algorithm involves symbolic computations since $c_2=0$.
\begin{table}[!h]
\caption{\label{tab1}The CPU times of \textbf{DETGTRI}, \textbf{Algorithm 2},	\textbf{Algorithm 3} and	\textbf{MATLAB} function($det()$) for Example 3.3}
\centering
\begin{tabular}{ccccc}
\hline
$n$ & \textbf{DETGTRI} &	\textbf{Algorithm 2}&	\textbf{Algorithm 3}&\textbf{MATLAB} ($det()$)\\
& CPU time(s) & CPU time(s)&CPU time(s) & CPU time(s) \\
\hline
10000& 0.782&	0.329& 0.078&	55.191
\\
20000&1.516&	0.672& 0.172&	342.727
\\
30000& 2.188&	1.063& 0.485&	1636.835
\\
40000&3.109&	1.297
& 0.500& --
\\
50000&3.906&	1.937
& 0.640& --
\\
100000&7.437&	4.218
& 0.796&
--\\
\hline
\end{tabular}
\end{table}
\\
Table 1 shows that the \textbf{Algorithm 3} is superior comparing with the  \textbf{DETGTRI} algorithm. The MATLAB function $det()$ has the largest CPU time between all algorithms.\\ 

\textbf{Example 3.4.} Consider the matrix $T_n$ given by:\\
\hspace*{3cm} $T_n=(t_{ij})\left[
\begin{array}{ccccccc}
    1 & 1 & 0 & \cdots & \cdots &\cdots & 0\\
   n-1 & 1 & 2 &\ddots &  & & \vdots \\
   0 & n-2 & 1 & 3& \ddots &  & \vdots\\
   \vdots &  &\ddots &\ddots &\ddots & & 0\\
   \vdots & & &&2 &1 & n-1\\ 
   0& \cdots & \cdots &\cdots & 0 & 1 & 1\\
 \end{array}
\right]_n$\\
\\
Consider $det(T_n)$. The \textbf{DETGTRI} algorithm gives:\\
$$c_k=\left\{
  \begin{array}{lr} 
      \qquad k & \text{if}\quad k\quad \text{is odd} \\
      -(n-k) & \text{if}\quad k\quad \text{is even}
      \end{array}
      \right.$$\\
Therefore,
$$det(T_n)=\prod_{r=1}^n c_r=\left\{
  \begin{array}{lr} 
      \qquad 0 & \text{if}\quad n\quad \text{is even} \\
      \frac{(-1)^{\frac{n-1}{2}}n!}{2^{n-1}}\binom{n-1}{\frac{n-1}{2}} & \text{if}\quad n\quad \text{is odd}
      \end{array}
      \right.$$\\
on simplification. Note that $det(T_n)=0$ when $n$ is even although $c_i\ne 0$ for $i=1, 2, \cdots, n-1$. This is because $c_n=0$.\\
In Table 2, we list some numerical results for  \textbf{DETGTRI} algorithm, \textbf{Algorithm 2} and	\textbf{Algorithm 3}. The superiority of \textbf{Algorithm 3} is obvious in Fig. 1. \\

\begin{table}[!h]
\caption{\label{tab2}Comparing \textbf{DETGTRI} algorithm,  \textbf{Algorithm 2} and 	\textbf{Algorithm 3} for Example 3.4}
\centering
\begin{tabular}{cccc}
\hline
$n$ & \textbf{DETGTRI} &	\textbf{Algorithm 2}&	\textbf{Algorithm 3}\\
& CPU time(s) & CPU time(s)& CPU time(s)\\
\hline
1000&0.109&	0.063 & 0.047\\
1500&0.140&	0.094 & 0.062\\
2000&0.156&	0.105&	0.078\\
2500&0.172&	0.125&	0.092\\
3000&0.250&	0.152&	0.128 \\
\hline
\end{tabular}
\end{table}
\begin{figure}[H]
\centering
\includegraphics[width=10cm, height=7cm] {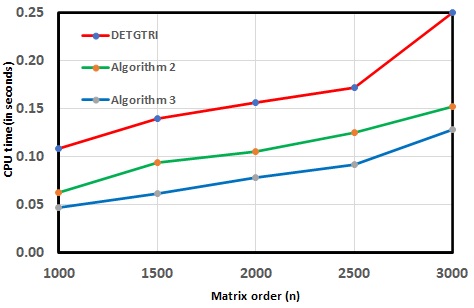}\qquad 
\caption{Efficiency of the \textbf{DETGTRI} algorithm, \textbf{Algorithm 2} and	\textbf{Algorithm 3}}
\label{Fig.1}
\end{figure}

\textbf{Example 3.5.} Consider the matrix $T_n$ given by:\\
\hspace*{3cm} $T_n=(t_{ij})\left[
\begin{array}{ccccccc}
    1 & 1 & 0 & \cdots & \cdots &\cdots & 0\\
   2 & 2 & 1 &\ddots &  & & \vdots \\
   0 & 2 & 2 & 1& \ddots &  & \vdots\\
   \vdots &  &\ddots &\ddots &\ddots & & 0\\
   \vdots & & &&2 &2 & 1\\ 
   0& \cdots & \cdots &\cdots & 0 & 2 & 1\\
 \end{array}
\right]_n$\\
\\
In this example, $c_2=0$. So, the \textbf{DETGTRI} algorithm contains symbolic computations. Table 3 shows the CPU times for the three algorithms. The \textbf{Algorithm 3} has CPU time less than the other two algorithms.
\begin{table}[!h]
\caption{\label{tab3}Comparing \textbf{DETGTRI} algorithm,  \textbf{Algorithm 2} and 	\textbf{Algorithm 3} for Example 3.5}
\centering
\begin{tabular}{cccc}
\hline
$n$ & \textbf{DETGTRI} &	\textbf{Algorithm 2}&	\textbf{Algorithm 3}\\
& CPU time(s) & CPU time(s)& CPU time(s)\\
\hline
1000&1.3440&	0.0437&	0.0031
\\
1500&1.6250&	0.0547&	0.0125
\\
2000&1.8600	&0.0626	&0.0140
\\
2500&2.4690	&0.0688&	0.0186
\\
3000&2.7650	&0.0985	&0.0265
 \\
\hline
\end{tabular}
\end{table}
Fig. 2 displays the logarithm of the CPU times multiplied by 1000 versus the matrix order $n$. Based on this figure, the \textbf{Algorithm 2} has least CPU times between all three algorithms.
\begin{figure}[H]
\centering
\includegraphics[width=10cm, height=7cm] {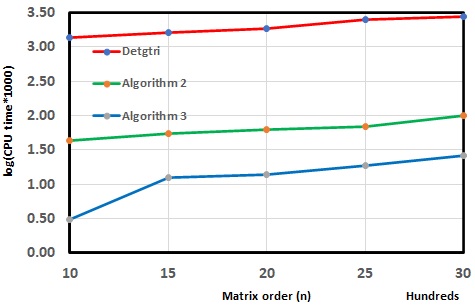}\qquad 
\caption{Efficiency of the \textbf{DETGTRI} algorithm, \textbf{Algorithm 2} and	\textbf{Algorithm 3}}
\label{Fig.2}
\end{figure}
\section{Conclusion}
In this paper, a hybrid numerical algorithm (Algorithm 3) has been derived for evaluating general n-th order tridiagonal determinants in linear time. The algorithm avoids all symbolic computations. The results show how effective the hybrid numerical algorithm is.

\bibliography{mybibfile}

\end{document}